\documentclass[a4paper,12pt]{article}
\usepackage{amsmath}
\usepackage{amsthm}
\usepackage{amssymb}
\usepackage{amsfonts}
\usepackage{enumerate}
\usepackage[all]{xy}

\normalsize

\newtheoremstyle{theoremstyle}
  {10pt}      
  {5pt}       
  {\itshape}  
  {}          
  {\bfseries} 
  {:}         
  {.5em}      
  {}          

\newtheoremstyle{examplestyle}
  {10pt}      
  {5pt}       
  {}          
  {}          
  {\bfseries} 
  {:}         
  {.5em}      
  {}          

\theoremstyle{theoremstyle}
\newtheorem{theorem}{Theorem}[section]
\newtheorem{lemma}[theorem]{Lemma}
\newtheorem{proposition}[theorem]{Proposition}
\newtheorem*{proposition*}{Proposition}
\newtheorem{corollary}[theorem]{Corollary}
\newtheorem*{corollary*}{Corollary}

\theoremstyle{examplestyle}
\newtheorem{example}[theorem]{Example}
\newtheorem{definition}[theorem]{Definition}
\newtheorem{remark}[theorem]{Remark}
\newtheorem{remark*}{Remark}

\newcommand{\comment}[1]{}

\newcommand{\rays}{{\Delta(1)}}
\newcommand{\mrays}{$\Delta(1)$}
\newcommand{\weildivisors}{{\mathbb{Z}^\rays}}
\newcommand{\mweildivisors}{$\weildivisors$}
\newcommand{\sh}[1]{\mathcal{#1}}
\newcommand{\msh}[1]{$\sh{#1}$}
\newcommand{\orb}[1]{\operatorname{orb}(#1)}
\newcommand{\fitt}[1]{\operatorname{Fitt}_2(#1)}
\newcommand{\PGL}{\operatorname{PGL}}
\newcommand{\GL}{\operatorname{GL}}
\newcommand{\gr}{\operatorname{Gr}}

\title
      {Moduli for Equivariant Vector Bundles of Rank Two on Smooth Toric Surfaces}
\author{Markus Perling\\ \footnotesize Department of Mathematics\\[-9pt]
        \footnotesize University of Kaiserslautern, Germany\\[-9pt] 
        \footnotesize \it email: \tt perling@mathematik.uni-kl.de}
\date{February 2002}

\begin{document}

\maketitle

\begin{abstract}
We give a complete classification of equivariant vector bundles of rank two over smooth
complete toric surfaces and construct moduli spaces of such bundles. This note is a direct
continuation of an earlier note where we developed a general description of equivariant
sheaves on toric varieties. Here we give a first application of that description.
\end{abstract}

\section{Introduction}

This note is a direct continuation of \cite{perling1} where, based on earlier work of
Klyachko (\cite{Kly90}, \cite{Kly91}), we have developed a formalism to describe equivariant
sheaves on toric varieties in terms of families of vector spaces and filtrations.
In \cite{perling1} we also started to consider free resolutions of fine-graded modules
over the homogeneous coordinate ring of a toric variety.
In fact, filtrations and free resolutions of fine-graded
modules over the homogeneous coordinate ring are in some sense dual notions. In this note
we want to give first examples of such a duality for the case of equivariant vector bundles
of rank
2 on smooth complete toric surfaces, which is the first nontrivial case one can consider.
Given a smooth complete toric surface $X$ which is determined by a fan $\Delta$, we
generalize a result of Kaneyama (\cite{Kan1}) and we will show
that every vector bundle \msh{E} of rank 2 can be realized as the cokernel in a
{\it generalized equivariant Euler type short exact sequence} (see Theorem
\ref{evenmoreeulersequences}):
\begin{equation*}
0 \longrightarrow \sh{O}^{s - 2} \overset{A}{\longrightarrow}
\bigoplus_{i = 1}^s  \sh{O}(\sum_{\rho \in \Pi_i} i^\rho \cdot D_\rho) \longrightarrow
\sh{E} \longrightarrow 0.
\end{equation*}
The precise shape of this sequence is determined by certain combinatorial data associated to
\msh{E} which will be expressed in terms of a {\it partition} $\{\Pi_i\}_{1 \leq i \leq s}$
of the set of rays in $\Delta$. Fixing this combinatorial data allows us to consider families
of equivariant torsion free sheaves of rank 2 on $X$ by varying matrices $A$ as in the
sequence above.
On the other hand, the above resolution of \msh{E} is derived from its representation in
terms of filtrations of a 2-dimensional vector space and indeed there is an immediate
relation between the space of variations of $A$ and the configuration space of filtrations
describing \msh{E}.

{\bf Plan of the paper:}
After recalling in Sections \ref{prereq} and \ref{eqsheaves} preliminaries on toric
geometry and equivariant sheaves, we briefly review in Section \ref{families} general
results on equivariant sheaves from \cite{perling1}. In Sections \ref{resolutions} and
\ref{partitionsresolutions}, we will use this description to construct and analyze
resolutions for general equivariant locally free sheaves of rank 2.
Section \ref{duality} is devoted to duality of configurations of points in the GIT setting
which we will use in Section
\ref{moduli} to give a GIT-classification of equivariant vector bundles of rank
2 on smooth complete toric surfaces.

{\bf Acknowledgements:} I want to express my gratitude to Prof. G. Trautmann for many discussions
and valuable advice.

\section{Toric Prerequisites}
\label{prereq}

The prerequisites for toric geometry are the same as in \cite{perling1} and we only want to
recall briefly some notions and to fix notation for the rest of this work. For the general
theory of toric varieties and standard notation we refer to the textbooks \cite{Oda} and
\cite{Fulton}.

\begin{itemize}
\setlength{\itemsep}{-5pt}
\item All algebraic varieties will be defined over a fixed algebraically closed field $k$,
\item $T \cong (k^*)^2$ denotes the 2-dimensional algebraic torus over $k$,
$M \cong \mathbb{Z}^2$ its character group, and $N$ the $\mathbb{Z}$-module dual to $M$;
the canonical pairing between $M$ and $N$ is denoted $\langle . , . \rangle$,
\item $X$ always denotes a smooth complete toric surface, described by a fan $\Delta$
\item cones in $\Delta$ are denoted by small Greek letters $\rho$, $\sigma$, $\tau$, etc.,
the natural order among cones is denoted by $\tau < \sigma$, \\
$\Delta(i) := \{\sigma \in \Delta \mid \dim \sigma = i\}$ the set of
all cones of fixed dimension $i$,
\item elements of $\Delta(1)$ are called {\it rays}, \\
the torus invariant Cartier divisor associated to the ray $\rho \in \rays$ is denoted
$D_\rho$, \\
$n(\rho) \in N$ denotes the primitive lattice element spanning the ray $\rho$,
\item $\check{\sigma} := \{m \in M_\mathbb{R} \mid \langle m, n \rangle \geq 0 \text{ for all
$n \in \sigma$}\}$ is the cone {\it dual} to $\sigma$, \\
$\sigma_M := \check{\sigma} \cap M$ denotes the subsemigroup of $M$ associated to $\sigma$,
\item the affine $T$-invariant subset associated to $\sigma \in \Delta$ is denoted
$U_\sigma \cong \operatorname{spec}(k[\sigma_M])$.
\end{itemize}

We will use {\it quotient presentations} $\hat{X} \overset{\pi}{\longrightarrow} X$,
due to Cox (\cite{Cox}). We consider the affine space $k^\rays$ as a toric variety
on which the torus $\hat{T} = (k^*)^\rays$ acts and denote $S = k[x_\rho \mid \rho \in
\rays]$ the ring of regular functions over $k^\rays$. $\hat{X}$ is the quasi-affine toric
subvariety of $k^\rays$ whose complement in $k^\rays$ is described by the $\hat{T}$-invariant
ideal $B$ in $S$, the so-called the {\it irrelevant ideal}. The codimension of $V(B)$
in $k^\rays$ is at least two, so that regular functions over $\hat{X}$ extend to regular
functions over $k^\rays$.  A theorem of Cox states that there is a morphism
$\hat{X} \overset{\pi}{\longrightarrow} X$
which is a categorical quotient with respect to the action of a diagonalizable subgroup
$G$ of $\hat{T}$ (see \cite{Cox} for details). The actions of $\hat{T}$ and $G$ both
induce gradings on $S$ with respect to their character groups $X(\hat{T})$ and $X(G)$
which are isomorphic to \mweildivisors\ and to the Chow group $A_{n - 1}(X)$, respectively.
Using this, one can consider homogeneous elements of $S$ with respect to the
$A_{n - 1}(X)$-grading
as global functions over $X$, and therefore $S$ is called the {\it homogeneous coordinate
ring of $X$}. Note that by the surjection $\weildivisors \twoheadrightarrow A_{n - 1}(X)$,
$\weildivisors$-homogeneneity implies $A_{n - 1}(X)$-homogeneity.
We will call the \mweildivisors-grading of $S$ the {\it fine grading} of $S$.

In our case of interest, where $X$ is a complete surface, we will make use of the following
explicit description of $B$.
Let $\rays = \{\rho_0, \dots, \rho_{m - 1}\}$ be the set of rays counted clockwise by the
cyclic group $\mathbb{Z}_m$. Then each 2-dimensional cone $\sigma \in \Delta$ is spanned
by rays $\rho_k$ and $\rho_{k + 1}$ for $k \in \mathbb{Z}_m$ and $B$ is generated by elements
$x^{\hat{\sigma}}$, i.e. $B = \langle x^{\hat{\sigma}} \mid \sigma \in \Delta(2) \rangle$,
where $x^{\hat{\sigma}} := \prod_{i = k + 2}^{k - 1}x_{\rho_i}$.

\section{Equivariant Sheaves}
\label{eqsheaves}

Let $G$ be an algebraic group which acts on a variety $X$. Following \cite{GIT}, \S 3, we call
a sheaf \msh{E} on $X$
{\it equivariant} or {\it $G$-linearized} if for each $g \in G$ considered as automorphism
of $X$, there is an isomorphism $\Phi_g : g^* \sh{E} \cong \sh{E}$ such that the following
diagram commutes for all $g, g' \in G$:

$$\xymatrix{
g^*{g'}^* \sh{E} \ar[rrd]_{g^*\Phi_{g'}} & & (g'g)^* \sh{E} \ar[ll]_\cong \ar@{-->}[d] \ar[rr]^{\Phi_{g'g}} & &  \sh{E} \\
& & {g}^* \sh{E} \ar[urr]_{\Phi_g}
}$$ 
We consider the case where $G$ is the algebraic torus $T$ and $X$
is a smooth complete toric surface over $k$.

One way to describe torus equivariant sheaves over $X$ is to consider
{\it fine-graded} modules over the homogeneous coordinate ring $S$ of $X$.
The sheafification operation which associates to each $A_{n - 1}(X)$-graded
$S$-module $F$ a quasicoherent sheaf $\tilde{F}$ over $X$, maps fine-graded
$S$-modules to equivariant sheaves over $X$.
We will use the following facts (see \cite{perling1} for detailed references):

\begin{theorem}
The map $F \mapsto \tilde{F}$ is a covariant additive exact functor from the category
of finitely generated fine-graded $S$-modules to the category of coherent equivariant
$\sh{O}_X$-modules. Moreover, every coherent equivariant sheaf is of the form $\tilde{F}$
for some finitely generated, fine-graded $S$-module $F$.
\end{theorem}

In particular, fine-graded free $S$-modules of rank 1, which are of the form
$S(\underline{n})$, where $\underline{n} = (n_\rho \mid \rho \in \rays) \in \weildivisors$
and $S(\underline{n})$ denotes the shift of degree by $\underline{n}$,
precisely correspond to equivariant line bundles $\sh{O}(\underline{n}) :=
\sh{O}(\sum_{\rho \in \rays} n_\rho \cdot D_\rho)$ over $X$.
Then an equivariant homomorphism of vector bundles
\begin{equation*}
\bigoplus_{j = 1}^m \sh{O}(\underline{m}_j) \overset{\tilde{A}}{\longrightarrow}
\bigoplus_{i = 1}^n \sh{O}(\underline{n}_i)
\end{equation*}
can be described by an $n \times m$-matrix $A = (a_{ij})$ whose entries are monomials
$a_{ij} = \alpha_{ij}x^{\underline{m}_j - \underline{n}_i}$, where $\alpha_{ij} \in k$
and $\alpha_{ij} = 0$ whenever
$\underline{m}_j - \underline{n}_i \notin \mathbb{N}^\rays$.

\section{$\Delta$-families and Filtrations}
\label{families}

We recall from \cite{perling1} some results on equivariant sheaves on toric varieties.
For each $\sigma \in \Delta$ we can define a preorder $\leq_\sigma$ on $M$ by setting
$m \leq_\sigma m'$ for $m, m' \in M$, iff
$m' - m \in \sigma_M$. In case that $m \leq_\sigma m'$ but not $m' \leq_\sigma m$
we write $m <_\sigma m'$ .

\begin{definition}
For a fixed $\sigma$ let $\{E^\sigma_m\}_{m \in M}$ be a family of $k$-vector spaces.
Let to each relation
$m \leq_\sigma m'$ a vector space homomorphism $\chi^{\sigma}_{m, m'} : E^\sigma_m
\longrightarrow E^\sigma_{m'}$ be given such that $\chi^\sigma_{m, m} = 1$ for
all $m \in M$ and $\chi^\sigma_{m, m''} = \chi^\sigma_{m', m''} \circ \chi^\sigma_{m, m'}$
for each triple $m \leq_\sigma m' \leq_\sigma m''$. We denote such a family by
$\hat{E}^\sigma$ and call it a {\it $\sigma$-family}.
\end{definition}

Any quasicoherent equivariant sheaf \msh{E} over $U_\sigma$ gives rise to a $\sigma$-family
by its isotypical decomposition which is induced by the action of $T$ on
$\Gamma(U_\sigma, \sh{E})$:
\begin{equation*}
\Gamma(U_\sigma, \sh{E}) = \bigoplus_{m \in M} \Gamma(U_\sigma, \sh{E})_m.
\end{equation*}
One just sets $E^\sigma_m = \Gamma(U_\sigma, \sh{E})_m$ and $\chi^\sigma_{m, m'}$
the monomial multiplication by $\chi(m' - m)$ which is canonically induced by the
ring structure. By \cite{perling1}, Prop. 5., the category of equivariant quasicoherent
sheaves over $U_\sigma$ is equivalent to the category of $\sigma$-families. This way, there
exists a natural process of gluing $\sigma$-families which mirrors the gluing of
quasicoherent sheaves. A collection $\{\hat{E}^\sigma\}_{\sigma \in \Delta}$ of
$\sigma$-families which glues is called a {\it $\Delta$-family} and denoted
$\hat{E}^\Delta$ (for details see \cite{perling1}).
From the functorial properties of gluing it follows that there
exists an equivalence of categories between $\Delta$-families and equivariant
quasicoherent sheaves over $X$.

Going one step further, given a $\Delta$-family $\hat{E}^\Delta$, for any
$\sigma \in \Delta$ we can consider the $\sigma$-family $\hat{E}^\sigma$ as a directed
system, and we  can form the direct limit $\underset{\rightarrow}{\lim}\hat{E}^\sigma =:
\mathbf{E}^\sigma$. Gluing then corresponds
to vector space isomorphisms $\mathbf{E}^\sigma \overset{\cong}{\longrightarrow}
\mathbf{E}^\tau$ whenever $\tau < \sigma$. In particular, for each $\sigma \in \Delta$
there is an isomorphism $\mathbf{E}^\sigma \overset{\cong}{\longrightarrow} \mathbf{E}^0$,
where $0$ is the minimal cone in $\Delta$. If $\hat{E}^\Delta$ defines a coherent sheaf
\msh{E}, then the $\mathbf{E}^\sigma$ are finite-dimensional of  dimension equal to
the rank of \msh{E}. If we restrict \msh{E} to the dense torus $T$, then $\sh{E}\vert_T$
is locally free and we can think of $\mathbf{E}^0$ as a general fibre of the associated
vector bundle over $T$.

Let us assume that $\hat{E}^\Delta$ represents a torsion free, coherent equivariant sheaf
over $X$. In that case, all the morphisms $\chi^\sigma_{m, m'}$ are injective and the
$E^\sigma_m$ become subvector spaces of $\mathbf{E}^\sigma$ such that morphisms
$\chi^\sigma_{m, m'}$ transform to inclusions $E^\sigma_m \subset E^\sigma_{m'} \subset
\mathbf{E}^\sigma$ whenever $m \leq_\sigma m'$.  This allows us to identify all the vector
spaces $E^\sigma_m$ as subvector spaces of the limit vector space $\mathbf{E}^0$, and the
category of torsion free, coherent equivariant sheaves over $X$ is equivalent to the
category of {\it embedded} $\Delta$-families. This category can be characterized by the
notion of {\it multifiltered} vector spaces, for whose somewhat lengthy definition we
refer to \cite{perling1}, definition 5.20. In this paper we will only need explicitly
the special case of reflexive sheaves. In that case, for $\sigma \in \Delta$ we can
reconstruct a $\sigma$-family $\hat{E}^\sigma$
from the $\rho$-families associated to the 1-dimensional cones $\rho$ of $\sigma$, via
$E^\sigma_m = \bigcap_{\rho \in \sigma(1)} E^\rho_m$.
In \cite{perling1} it is shown that a $\rho$-family for $\rho \in \rays$ is equivalent
to a full filtration of the vector space $\mathbf{E}^0$:
\begin{equation*}
\dots \subset E^\rho(i) \subset E^\rho(i + 1) \subset \dots \subset \mathbf{E}^0
\end{equation*}
where 'full' means that there exist numbers $i_1 < i_2$ such that $E^\rho(i) = 0$ for
$i < i_1$ and $E^\rho(i) = \mathbf{E}^0$ for $i > i_2$.
There is an identification of $E^\rho_m$ with $E^\rho(\langle m, n(\rho) \rangle)$,
where $n(\rho)$ is the primitive lattice vector of the ray $\rho$.

By the functoriality of direct limits an equivariant homomorphism of torsion free sheaves
$\sh{E} \longrightarrow \sh{E}'$ is equivalent to a vector space
homomorphism $\mathbf{E}^0 \longrightarrow (\mathbf{E}')^0$ compatible with the
multifiltrations of $\mathbf{E}^0$ and $(\mathbf{E}')^0$.
We summarize all this in:

\begin{theorem}[\cite{perling1}, Theorem 5.22]
let $X$ be a toric variety. The category of equivariant reflexive sheaves over $X$ is equivalent to
the category of vector spaces with full filtrations associated to each ray in \mrays\ and
whose morphisms are vector space homomorphisms which are compatible with the associated
$\Delta$-families.
\end{theorem}
In particular, if $X$ is a smooth toric surface, then all reflexive sheaves are locally
free, and thus the category of filtrations is equivalent to that of equivariant
locally free sheaves.

\section{Partitions and Resolutions of Equivariant Bundles of Rank Two}
\label{partitions}
\label{resolutions}

In this section we want to construct resolutions for general equivariant
vector bundles of rank 2 on toric surfaces. Let \msh{E} be such a bundle on $X$ given
by filtrations of a 2-dimensional vector space $\mathbf{E}^0$. Any such filtration
$E^\rho(i)$ can be described as follows. There are two integers $i^\rho_1 \leq i^\rho_2$
such that
\begin{equation*}
E^\rho(i) =
\begin{cases}
0 & \text{\quad for $i < i^\rho_1$},\\
E^\rho & \text{\quad for $i^\rho_1 \leq i < i^\rho_2$},\\
\mathbf{E}^0 & \text{\quad for $i \geq i^\rho_2$},
\end{cases}
\end{equation*}
where $E^\rho$ is a 1-dimensional subvector space of $\mathbf{E}^0$. Thus in the case that
$i^\rho_1 < i^\rho_2$ a filtration can irredundantly be described
by an ordered triple $(i^\rho_1, i^\rho_2, E^\rho)$. If $i^\rho_1 = i^\rho_2$,
the filtration is degenerate in the sense that at $i^\rho_1$ the dimension jumps by two
and there occurs no $E^\rho$.

Twisting \msh{E} by an equivariant line bundle $\sh{O}(\underline{n})$ for some
$\underline{n} = (n_\rho) \in \weildivisors$ has the
effect that the numbers $i^\rho_1, i^\rho_2$ are shifted to $i^\rho_1 + n_\rho,
i^\rho_2 + n_\rho$ for all $\rho$ (see \cite{perling1}, \S 6). In our further
considerations such twists do not play any role and we may assume for
simplicity that $i^\rho_1 = -i^\rho$ and $i^\rho_2 = 0$ for nonnegative
integers $i^\rho$.  In \cite{perling1} we have obtained the following result
which we are going to generalize:

\begin{theorem}[\cite{perling1}, Theorem 6.1]
\label{eulersequences}
Let $X$ be a complete smooth toric surface and let \msh{E} be an equivariant rank 2-vector
bundle over $X$ determined by filtrations $\{(-i^\rho, 0, E^\rho)\}_{\rho \in \rays}$.
If $i^\rho > 0$ for
all $\rho \in  \rays$ and $E^\rho \neq E^{\rho'}$ whenever there is a
$\sigma \in \Delta(2)$ such that $\rho, \rho' \in \sigma(1)$, then there exists an
Euler type equivariant short exact sequence:
\begin{equation*}
0 \longrightarrow \sh{O}^{n - 2} \overset{A}{\longrightarrow} \bigoplus_{\rho \in \rays}
\sh{O}(i^\rho \cdot D_\rho) \longrightarrow \sh{E} \longrightarrow 0
\end{equation*}
where $n = \# \rays$ and $A = (\alpha_{\rho j} \cdot x_\rho^{i_\rho})$, $1 \leq j \leq n - 2$,
$\alpha_{\rho j} \in k$, an $n \times (n - 2)$-matrix of monomials.
\end{theorem}

This theorem yields resolutions for {\it generic} equivariant vector bundles \msh{E} of
rank 2. In order to obtain a complete classification and moduli, one has also to
consider the two ways how a set of filtrations may degenerate. Namely, the case where
$i^\rho = 0$ for some $\rho \in \rays$ and $E^\rho = E^{\rho'}$ for two adjacent
rays $\rho$ and $\rho'$

In order to do this, we consider {\it partitions} of subsets $\Pi \subset \rays$ as follows.
A partition of $\Pi$ is a collection $\mathcal{P} = \{\Pi_1, \dots, \Pi_s\}$ of disjoint
subsets of $\Pi$ such that $\Pi = \coprod_{i = 1}^s \Pi_i$. Among the partitions of
$\Pi$ we define a partial order $\leq$ as follows: given two partitions $\mathcal{P} =
\{\Pi_1, \dots, \Pi_s\}$
and $\mathcal{P}' = \{\Pi'_1, \dots, \Pi'_{s'}\}$ we say $\mathcal{P} \leq \mathcal{P}'$
iff $\mathcal{P}$ is a {\it refinement} of $\mathcal{P}'$, i.e. there exists a map
$\{1, \dots, s\} \longrightarrow \{1, \dots s'\}$, $i \mapsto j$ such that $\Pi_i \subset
\Pi'_j$. We call the map given by a refinement
\begin{equation*}
\pi: \mathcal{P} \longrightarrow \mathcal{P}', \quad \pi(\Pi_i) = \Pi'_j,
\end{equation*}
{\it projection map},
and any map $s : \mathcal{P}' \longrightarrow \mathcal{P}$ such that $\pi \circ s$
is the identity a {\it section} with respect to $\pi$. The partial order $\leq$ has
unique minimal and maximal elements, namely the partitions
$\{\{\rho\}\}_{\rho \in \Pi}$ and $\{\Pi\}$.

\paragraph{Partitions associated to an equivariant bundle of rank two:}
For an equivariant bundle \msh{E} of rank 2 we denote by $\Pi = \Pi(\sh{E}) \subset \rays$ the
subset of those $\rho$ for which $i^\rho > 0$.
We assume that the rays $\{\rho_0, \dots \rho_{m - 1}\} = \Pi$ are enumerated clockwise 
with respect to their circular order in the fan by the
cyclic group $\mathbb{Z}_m$. On $\Pi$ we can define a partition as follows.

\begin{definition}
Let $\Pi_1, \dots \Pi_s$ be the unique partition of $\Pi$ with the following properties:
\begin{enumerate}[(i)]
\item if $\rho_i, \rho_j \in \Pi_k$ then $E^{\rho_i} = E^{\rho_j}$,
\item if for some $1 \leq k < s$ $\rho_i \in \Pi_k$ and $\rho_j \in \Pi_{k + 1}$,
or if $\rho_i \in \Pi_s$ and $\rho_j \in \Pi_1$, then $E^{\rho_i} \neq E^{\rho_j}$,
\item if $\Pi_k$ contains $\rho_i$ and $\rho_{i + l}$ then it contains the interval $\rho_{i + 1},
\dots, \rho_{i + l - 1}$ or the interval $\rho_{i + l + 1}, \dots, \rho_{i - 1}$ or
both.
\end{enumerate}
One can think of this partition as the set of maximal intervals in the circularly
ordered set $\Pi$ on which the $E^\rho$ coincide.
We assume that the $\Pi_i$ are enumerated clockwise. We denote the partition so defined
$\mathcal{P}_\mathcal{E}$ and call it {\it the coarse partition of $\Pi$ with
respect to \msh{E}}.
\end{definition}

In the situation of Theorem \ref{eulersequences}, the coarse partition
$\mathcal{P}_\mathcal{E}$ of \mrays\ of a bundle \msh{E} is precisely
the partition $\{\{\rho\}\}_{\rho \in \rays}$. Using the definition of coarse partitions, 
we can extend Theorem \ref{eulersequences}
to the case of any equivariant bundle of rank 2 on $X$.

\begin{theorem}
\label{evenmoreeulersequences}
Let \msh{E} be an arbitrary equivariant vector bundle of rank 2 on a smooth complete
toric surface
$X$, defined by filtrations $\{(-i^\rho, 0, E_\rho)\}_{\rho \in \rays}$ of a two
dimensional vector space $\mathbf{E}^0$. Let $\Pi =
\{\rho \in \rays \mid i^\rho > 0\}$ and let
$\mathcal{P}_\mathcal{E} = \{\Pi_1, \dots \Pi_s\}$ be the coarse partition of
$\Pi$ with respect to \msh{E}.
If $s > 2$ then there exists a short exact sequence
\begin{equation*}
0 \longrightarrow \sh{O}^{s - 2} \overset{A}{\longrightarrow}
\bigoplus_{i = 1}^s  \sh{O}(\sum_{\rho \in \Pi_i} i^\rho \cdot D_\rho) \longrightarrow
\sh{E} \longrightarrow 0
\end{equation*}
where $A$ is a matrix of monomials whose exponents are determined by the partition
$\mathcal{P}_\mathcal{E}$. Moreover, the $(s - 2)$-minors $A^{i,i+1}$ of $A$, $1 \leq i < s$,
which consist of all
rows of $A$ except the $i$-th and the $(i+1)$-st, are of full rank.
If $s \leq 2$, then \msh{E} splits.
\end{theorem}

\begin{remark}
The proof explains the precise relationship between $A$ and the filtrations associated to
\msh{E}, see also Proposition \ref{desingprep1}.
\end{remark}

\begin{proof}
Let first $s \leq 2$. Then we can decompose the vector space $\mathbf{E}^0$ into a direct sum
$\mathbf{E}^0 = E_1 \oplus E_2$
and the filtrations decompose into direct sums of filtrations for $E_1$ and $E_2$,
respectively. Consequently, the associated bundle \msh{E} splits into a direct sum
of line bundles.

Now assume that $s > 2$.
Consider the Cox quotient presentation $\hat{X} \longrightarrow X$ and let
$\{n_\rho\}_{\rho \in \rays}$ be the standard basis of the lattice $\hat{N} \cong
\weildivisors$. Let $\{i_1, \dots, i_{s - 2}\} \subset \{1, \dots, s\}$, then we set
\begin{equation*}
x^{\hat{\Pi}_i} := \prod_{\rho \in \Pi_{i}} x^{i^\rho}_\rho, \text{\quad and \quad}
x^{\hat{\Pi}_{i_1 \dots i_{s - 2}}} := \prod_{k = 1}^{s - 2} x^{\hat{\Pi}_{i_k}}.
\end{equation*}
We can define a morphism of fine-graded free $S$-modules
\begin{equation*}
0 \longrightarrow S^{s-2} \overset{A}{\longrightarrow} \bigoplus_{i = 1}^s
S(\sum_{\rho \in \Pi_i} i^\rho \cdot n_\rho)
\end{equation*}
which is an $s \times (s - 2)$-matrix $A$ with monomial entries:
\begin{equation*}
A = (\alpha_{ij} \cdot x^{\hat{\Pi}_i})
\end{equation*}
where $i$ runs from $1$ to $s$ and $j$ from $1$ to $s - 2$. We require that, for
$i = 1, \dots, s - 1$, the $(s - 2)$-minors $A^{i, i+1}$ of $A$ which
consist of all rows of $A$ except the $i$-th and the $(i+1)$-st, have full rank over $S$.
After applying the sheafification functor $\tilde{\ }$
to this sequence we obtain a short exact sequence of sheaves
\begin{equation*}
0 \longrightarrow \sh{O}^{s - 2} \overset{A}{\longrightarrow} \bigoplus_{i = 1}^s
\sh{O}(\sum_{\rho \in \Pi_i} i^\rho \cdot D_\rho) \longrightarrow \sh{Q} \longrightarrow 0
\end{equation*}
where by abuse of notion we write $A$ instead of $\tilde{A}$.

The matrix $A$ defines an equivariant injective morphism of coherent sheaves,
but it is not necessarily an injective vector bundle homomorphism. This is the case if
and only if the rank of $A(x)$
equals $s - 2$ at all points $x \in X$. This in turn means that $A$ is an inclusion of
vector bundles if and only if there exists a $k > 0$ such that $B^k \subset \fitt{A}$,
where $B$ is the irrelevant ideal associated to the quotient presentation $\hat{X}
\longrightarrow X$. If this is
the case, then the cokernel \msh{Q} is a vector bundle as well.

Let $\{i_1, \dots, i_{s - 2}\} \subset \{1, \dots s\}$ and let $A^{i_1 \dots i_{s - 2}}$
be the $(s - 2) \times (s - 2)$-minor of $A$ which contains the rows corresponding
to $\{i_1, \dots i_{s - 2}\}$. Moreover, let
\begin{equation*}
A' := (\alpha_{i, j})
\end{equation*}
be the matrix of coefficients of $A$ and $(A')^{i_1 \dots i_{s - 2}}$ the according minor.
The second Fitting ideal $\fitt{A}$ of $A$ is generated by the determinants of all the
$A^{i_1 \dots i_{s - 2}}$:
\begin{equation*}
\fitt{A} = \langle \det A^{i_1 \dots i_{s - 2}} \rangle =
\langle \det (A')^{i_1 \dots i_{s - 2}} \cdot x^{\hat{\Pi}_{i_1 \dots i_{s - 2}}} \rangle
\end{equation*}
Thus $\fitt{A}$ is a monomial ideal generated by the $x^{\hat{\Pi}_{i_1 \dots i_{s - 2}}}$.
To show that $B^k \subset \fitt{A}$ for some $k > 0$, it suffices to show that for each
generator $x^{\hat{\sigma}}$ of $B$ there exists a generator of $\fitt{A}$ which divides
some power of $x^{\hat{\sigma}}$. Without loss of generality we may assume that $i^\rho = 1$
for all $\rho \in \Pi$. Then the problem is equivalent to the question whether a given
$x^{\hat{\sigma}}$ with $\sigma(1) = \{\rho_k, \rho_{k + 1}\}$ is divided by some
$x^{\hat{\Pi}_{i_1 \dots i_{s - 2}}}$ which in turn is equivalent to finding
$i_1, \dots i_{s - 2}$ such
that $\Pi_{i_1} \cup \dots \cup \Pi_{i_{s - 2}}$ is contained in the interval $\{\rho_{k + 2}
\dots \rho_{k - 1}\}$ (with indices modulo $m$). But because the complement of this interval
is $\{\rho_k, \rho_{k + 1}\}$, this complement intersects at most two of the intervals
$\Pi_i$, say, after renumbering, $\Pi_{s - 1}$ and $\Pi_s$. Thus we choose
$(i_1, \dots i_{s - 2}) = (1, \dots, s - 2)$ which is nonempty because $s > 2$.
Moreover, $\Pi_1 \cup \dots \cup \Pi_{s - 2}$, is contained in $\{\rho_{k + 2}, \dots,
\rho_{k - 1}\}$, and so $x^{\hat{\Pi}_{i_1, \dots, i_{s - 2}}}$ divides $x^{\hat{\sigma}}$.
 Hence, chosing a matrix $A$ as
above ensures that the quotient \msh{Q} of $A$ is locally free.

Now we have to show that any \msh{E} with associated coarse partition
$\mathcal{P}_\mathcal{E}$ can be resolved this way. We do this by explicitly
writing down the filtrations for $\sh{O}^{s - 2}$ and
$\bigoplus_{i = 1}^s \sh{O}(\sum_{\rho \in \Pi_i}i^\rho \cdot D_\rho)$ and by constructing
with their help a homomorphism of the associated limit vector spaces which we will lift
to a morphism of locally free sheaves.
Denote $\mathbf{F}^0$ the $(s - 2)$-dimensional filtered $k$-vector space associated to the
vector bundle $\sh{O}^{s-2}$, and $\mathbf{G}^0$ the $s$-dimensional $k$-vector space
associated to $\bigoplus_{i = 1}^s \sh{O}(\sum_{\rho \in \Pi_i}i^\rho \cdot D_\rho)$.
We will identify $\mathbf{G}^0$ with $k^{\mathcal{P}_\sh{E}} \cong k^s$ and label its
standard basis $e_1, \dots, e_s$. The filtrations are:
\begin{equation*}
F^\rho(i) =
\begin{cases}
0 & \text{\quad for } i < 0 \\
F & \text{\quad otherwise}
\end{cases}
\text{\quad and }
G^\rho(i) = 
\begin{cases}
0 & \text{\quad for } i < -i^\rho \\
k \cdot e_i & \text{\quad for } -i^\rho \leq i < 0 \text{ and } \rho \in \Pi_i\\
\mathbf{G}^0 & \text{\quad otherwise}
\end{cases}
\end{equation*}
The matrix $A$ induces a vector space homomorphism from $\mathbf{F}^0$ to $\mathbf{G}^0$
which can be naturally identified with the matrix $A'$.
We can define filtrations for the quotient vector space
$\mathbf{E}^0 := \mathbf{G}^0 / \mathbf{F}^0$ simply by taking the quotient filtrations
\begin{equation*}
E^\rho(i) = G^\rho(i) / F^\rho(i)
\end{equation*}
with respect to $A'$. These filtrations are of the form
\begin{equation*}
E^\rho(i) = (-i^\rho, 0, k \cdot \overline{e}_j)
\end{equation*}
where $\rho \in \Pi_j$ and $\overline{e}_j$ is the image of $e_j$ in $\mathbf{E}^0$.
If we assume that $B^k \subset
\fitt{A}$ for some $k > 0$, these filtrations become in a natural way the filtrations
associated to the cokernel \msh{E} of $A$.
On the other hand, if we define a homomorphism from $\mathbf{G}^0$ to some 2-dimensional
$k$-vector space $\mathbf{E}^0$ by fixing the images $\overline{e}_j \neq 0$ of the basis
vectors $e_j$, $j = 1, \dots, s$, of $\mathbf{G}^0$, we immediately
obtain a homomorphism of filtered vector spaces whose kernel is a filtered vector space
$\mathbf{F}^0 \overset{A'}{\hookrightarrow} \mathbf{G}^0$.
The corresponding matrix $A$ with monomial entries then defines a sheaf homomorphism
$0 \longrightarrow \sh{O}^{s - 2} \overset{A}{\longrightarrow} \bigoplus_{i = 1}^s
\sh{O}(\sum_{\rho \in \Pi_i} i^\rho \cdot D_\rho)$. As we have seen before, the cokernel of
$A$ is locally free
if and only if $\det(A)^{i, i+1} \neq 0$ for $i = 1, \dots, s - 1$ and $\det(A)^{s,1} \neq 0$.
Now it is a lemma from linear algebra that $\overline{e}_i$ and $\overline{e}_{i + 1}$ are
linearly independent if and only if $\det (A')^{i, i + 1} \neq 0$.
\end{proof}

\section{More on Partitions and Resolutions}
\label{partitionsresolutions}

Let us fix numbers $i^\rho > 0$ for $\rho \in \Pi \subset \rays$ and a partition $\mathcal{P}$
of $\Pi$. In this section we consider short exact sequences of type
\begin{equation*}
0 \longrightarrow \sh{O}^{s - 2} \overset{A}{\longrightarrow}
\bigoplus_{i = 1}^s \sh{O}(\sum_{\rho \in \Pi_i} i^\rho \cdot D_\rho)
\longrightarrow \sh{E} \longrightarrow 0
\end{equation*}
where $A$ is given by a monomial matrix.
By Theorem \ref{evenmoreeulersequences}, there are conditions such that the cokernel
\msh{E} is a locally free sheaf whose associated coarse partition $\mathcal{P}_\mathcal{E}$
coincides with $\mathcal{P}$. In general, if $A$ is arbitrary and has just maximal rank, we
have the following as an immediate corollary from the constructions of the previous
section:

\begin{proposition}
\label{desingprep1}
Fix a set of numbers $I := \{i^\rho \geq 0\}_{\rho \in \rays}$, let $\Pi = \Pi_I = \{\rho \mid i^\rho > 0\}
\subset \rays$ and let $\mathcal{P} = \{\Pi_1, \dots, \Pi_s\}$, $s \geq 2$, be a partition
of $\Pi$. Let $\sh{E} = \sh{E}(I, \mathcal{P}, A)$ be a sheaf defined by a short exact sequence
\begin{equation*}
0 \longrightarrow \sh{O}^{s - 2} \overset{A}{\longrightarrow}
\bigoplus_{i = 1}^s \sh{O}(\sum_{\rho \in \Pi_i} i^\rho \cdot D_\rho)
\longrightarrow \sh{E} \longrightarrow 0.
\end{equation*}
Then \msh{E} is a torsion free sheaf of rank 2 over $X$ and
we can consider the short exact sequence induced on the limit vector spaces
\begin{equation*}
0 \longrightarrow k^{s - 2} \overset{A^0}{\longrightarrow} k^{\mathcal{P}}
\overset{\check{A}^0}{\longrightarrow} \mathbf{E}^0 \longrightarrow 0,
\end{equation*}
with $\check{A}^0 = (A_1, \dots, A_s)$ a $2 \times s$-matrix. If \msh{E} is locally free, then
$\mathcal{P} \leq \mathcal{P}_\mathcal{E}$ is a refinement of the coarse partition
associated to \msh{E}, and the
filtrations for \msh{E} are given by $\{(-i^\rho, 0, \langle A_i \rangle) \mid \rho \in \Pi_i\}_{\Pi_i \in
\mathcal{P}}$,
where $\langle A_i \rangle$ denotes the 1-dimensional subvector space of
$\mathbf{E}^0$ spanned by the $i$-th column of $\check{A}^0$.
\end{proposition}

\begin{remark}
\label{desingprep2}
Observe that for $A$ and \msh{E} as in Proposition \ref{desingprep1} and \msh{E} locally free,
and for $\{\Pi'_1, \dots \Pi'_{s'}\} = \mathcal{P}' \leq \mathcal{P}$ any refinement with corresponding
projection $\pi$, we can write the filtrations as $\{(-i^\rho, 0, \langle A_{\pi(i)} \rangle) \mid \rho \in
\Pi'_i\}_{\Pi'_i \in \mathcal{P}'}$.
\end{remark}

Using the fact that each torsion free sheaf \msh{E} embeds into its bidual, $0 \longrightarrow
\sh{E} \longrightarrow \sh{E}\check{\ }\check{\ }$, we can now completely describe torsion free
equivariant sheaves of rank 2 over $X$ without explicitly considering $\Delta$-families:

\begin{theorem}
\label{desing}
Let $\sh{E}' = \sh{E}'(I, \mathcal{P}', B)$ be a cokernel
\begin{equation*}
0 \longrightarrow \sh{O}^{s' - 2} \overset{B}{\longrightarrow}
\bigoplus_{i = 1}^{s'} \sh{O}(\sum_{\rho \in \Pi'_i} i^\rho \cdot D_\rho)
\longrightarrow \sh{E}' \longrightarrow 0.
\end{equation*}
and let $\check{B}^0$ be defined by
\begin{equation*}
0 \longrightarrow k^{s' - 2} \overset{B^0}{\longrightarrow} k^{\mathcal{P}'}
\overset{\check{B}^0}{\longrightarrow} k^2 \longrightarrow 0,
\end{equation*}
Let then \msh{E} be the bundle defined by the filtrations $\{-i^\rho, 0, E^\rho\}_{\rho \in \Pi}$
associated to $\check{B}^0$ by
\begin{equation*}
E^\rho = \begin{cases}
0 & \rho \notin \Pi \\
\langle B_i \rangle & \rho \in \Pi'_i
\end{cases}
\end{equation*}
Then $\sh{E} \cong \sh{E}'\check{\ }\check{\ }$, $\mathcal{P}' \leq \mathcal{P}_\mathcal{E}$, and
we have an exact diagram
\begin{equation*}
\label{reductiondiagram}
\xymatrix{
& &  & 0 \ar[d] \\
0 \ar[r] & \sh{O}^{s' - 2} \ar[r]^-B \ar[d] & \bigoplus_{i = 1}^{s'}
\sh{O}(\sum_{\rho \in \Pi'_i}i^\rho \cdot D_\rho) \ar[r]^-{\check{B}} \ar[d]^\pi & \sh{E}'
\ar[r] \ar[d] & 0 \\
0 \ar[r] & \sh{O}^{s - 2} \ar[r]^-A \ar[d] & \bigoplus_{i = 1}^s
\sh{O}(\sum_{\rho \in \Pi_i} i^\rho \cdot D_\rho)
\ar[r]^-{\check{A}} \ar[d] & \sh{E} \ar[r] \ar[d] & 0 \\
& 0 \ar[r] & \mathcal{C} \ar[r] \ar[d] & \mathcal{C} \ar[r] \ar[d] & 0 \\
& & 0 & 0
}
\end{equation*}
The cokernel sheaf \msh{C} is a skyscraper sheaf whose
support is contained in the set of 0-dimensional orbits of $X$. More precisely,
\begin{equation*}
\operatorname{supp}(\sh{C}) = \dot{\bigcup_{\sigma \in \delta(2)}} \orb{\sigma}
\end{equation*}
where $\delta(2) \subset \Delta(2)$ is the set of cones $\sigma \in \Delta(2)$ such
that for $\{\rho_i, \rho_{i + 1}\} = \sigma(1)$, it is true that $\rho_i \in \Pi'_i$,
$\rho_{i + 1} \in \Pi'_j$, for some $i \neq j$, and $\pi(\Pi'_i) = \pi(\Pi'_j)$.
\end{theorem}

\begin{proof}
Denote $\mathcal{P} := \mathcal{P}_\mathcal{E}$ and let
$\mathcal{P}' \leq \mathcal{P}$ be any refinement with projection
$\pi$. Using a section $t : \mathcal{P}' \longrightarrow \mathcal{P}$ we fix a choice of elements
in the preimage of $\pi$. We define a matrix $\check{A}^0 := (B_{t(1)}, \dots, B_{t(s)})$ and 
$\hat{\pi} : k^{\mathcal{P}'} \longrightarrow k^{\mathcal{P}}$ the morphism induced by $\pi$ over
$k$. This way we obtain in the category of $k$-vector spaces a commutative diagram
\begin{equation*}
\xymatrix{
0 \ar[r] & k^{s' - 2} \ar[r]^{B^0} \ar[d] & k^{\mathcal{P}'} \ar[r]^{\check{B}^0}
\ar[d]_{\tilde{\pi}} & \mathbf{E}^0 \ar[r] \ar[d]^{\operatorname{id}} & 0 \\
0 \ar[r] & k^{s - 2} \ar[r]^{A^0} & k^{\mathcal{P}} \ar[r]^{\check{A}^0} &
\mathbf{E}^0 \ar[r] & 0
}
\end{equation*}
where $\operatorname{id}^0$ is the identity homomorphism on  $\mathbf{E}^0$ and
$A^0$ the kernel homomorphism of $\check{A}^0$.

The morphisms in the left square of the diagram can immediately be lifted to morphisms of
locally free sheaves by considering them as matrices of coefficients of the entries of matrices
of monomials. So we obtain the diagram
\begin{equation*}
\xymatrix{
0 \ar[r] & \sh{O}^{s' - 2} \ar[r]^-B \ar[d] & \bigoplus_{i = 1}^{s'}
\sh{O}(\sum_{\rho \in \Pi'_i}i^\rho \cdot D_\rho) \ar[r]^-{\check{B}} \ar[d]^\pi & \sh{E}'
\ar[r] \ar[d] & 0 \\
0 \ar[r] & \sh{O}^{s - 2} \ar[r]^-A & \bigoplus_{i = 1}^s
\sh{O}(\sum_{\rho \in \Pi_i} i^\rho \cdot D_\rho)
\ar[r]^-{\check{A}} & \sh{E} \ar[r] & 0
}
\end{equation*}
where we interprete the matrices $\check{A}^0$ and $\check{B}^0$ as sheaf
homomorphisms. 
The injectivity of the homomorphism $\sh{E}' \longrightarrow \sh{E}$ is an immediate
consequence of the fact that after restriction to the open sets $U_\rho$, $\rho \in \rays$,
it induces the identity homomorphism. It follows that the
cokernel $\sh{C}$ is a skyscraper sheaf whose support must be contained in the set of
0-dimensional orbits of $X$, and its description is immediate.
\end{proof}

\section{Duality for Configuration Spaces of Points in Projective Spaces}
\label{duality}

Let $m < n$, $T_n \cong (k^*)^n$ the $n$-dimensional
algebraic torus, $\GL_m$ the group of automorphisms of $k^m$ and denote $G := \GL_m \times T_n$,
which is a reductive group. Denote $\mathbb{M}_{n,m}$ the space $\operatorname{Hom}_k(k^m,
k^n)$ of $n \times m$-matrices
over $k$ and let $G$ act on $\mathbb{M}_{n,m}$ by $(g, t) . A := t \circ A \circ g^{-1}$.

We first want to consider the actions of the two subgroups $\GL_m$ and $T_n$ of $G$
separately. Because the representations of $\GL_m$ and $T_n$ in $\GL(\mathbb{M}_{n,m})$
both contain the homotheties, their actions induce actions on the projective space
$\mathbb{P}\mathbb{M}_{n, m}$ and linearizations of the ample line bundle
$\sh{O}_{\mathbb{P}\mathbb{M}_{n, m}}(1)$, so that we are able to perform
GIT-quotients of $\mathbb{P}\mathbb{M}_{n, m}$ by $\GL_m$ and $T_n$, respectively.
$\GL_m$ acts from the right on the matrices $\mathbb{M}_{n,m}$, and the set of semistable
points in $\mathbb{P}\mathbb{M}_{n,m}$ is precisely the set of points represented by matrices
which have maximal rank $m$:
\begin{equation*}
\mathbb{P}\mathbb{M}_{n, m}^{ss}(\GL_m) = \{\langle A \rangle \mid \operatorname{rank}A = m \}
\end{equation*}
Furthermore, $\GL_m$ acts freely on this set, so that
$\mathbb{P}\mathbb{M}_{n, m}^{s}(\GL_m) = \mathbb{P}\mathbb{M}_{n, m}^{ss}(\GL_m)$,
and there exists the geometric quotient
\begin{equation*}
\mathbb{P}\mathbb{M}_{n,m}^{ss}(\GL_m) // \GL_m \cong \gr(m, n)
\end{equation*}
where $\gr(m, n)$ is the Grassmannian of $m$-dimensional linear subspaces of $k^n$.
Similarly, with help of the eigenspace decomposition of the left action of $T_n$
on $\mathbb{M}_{n,m}$, it easy to see that
\begin{equation*}
\mathbb{P}\mathbb{M}_{n, m}^{ss}(T_n) = \{\langle A \rangle \mid
\text{ no row of $A$ is zero }\}.
\end{equation*}
$T_n$ acts freely on this set, so that stable and semistable points of
$\mathbb{P}\mathbb{M}_{n, m}$ coincide and and we obtain a geometric quotient
\begin{equation*}
\mathbb{P}\mathbb{M}_{n,m}^{ss}(T_n) / T_n \cong (\mathbb{P}_{m - 1})^n
\end{equation*}
which is given by the map $\langle A \rangle \mapsto (\langle A_1 \rangle, \dots, \langle A_n
\rangle)$, where $A_i$ denote the row vectors of the matrix $A$.

The action of the group $G$ descends to actions of the groups $G / \GL_m \cong T_n$ on
$\gr(m, n)$ and $G / T_n \cong \GL_m$ on $(\mathbb{P}_{m - 1})^n$, respectively. Both
actions are textbook examples from GIT and there are the following criteria for
stability:
\begin{proposition}[\cite{GIT}, Proposition 4.3]
\label{tuplequot}
\ 

\begin{enumerate}
\item An $n$-tuple $(p_1, \dots, p_n)$ of points in $(\mathbb{P}_{m - 1})^n$ is
(semi-)stable with respect to the diagonal action of $\GL_m$
if and only if for every proper linear subspace $L$ of $k^m$
\begin{equation}
\label{tuplequotcondition}
\#\{i \mid p_i \in L \} < \frac{m}{n} \dim L
\end{equation}
(respectively $\leq$).
\item Consider the action of $\GL_n$ on $\gr(m, n)$.
Then a point $A \in \gr(m, n)$ is (semi-) stable with respect to this action if and only if,
for every proper linear subspace $L$ of $k^n$,
\begin{equation}
\label{grasstability}
\dim (A \cap L) < \frac{m}{n} \dim L
\end{equation}
(respectively $\leq$).
\end{enumerate}
\end{proposition}
We need to modify the second
statement only slightly for the case of the action of a maximal subtorus of $\GL_n$ on
$\gr(m, n)$:
\begin{corollary}
Consider the action of a maximal subtorus $T_n$ of $\GL_n$ on $\gr(m, n)$. 
Then a point $A \in \gr(m, n)$ is (semi-)stable with respect to this action if and only if,
for every proper linear subspace $L$ of $k^n$ which is spanned by eigenspaces of the
action of $T_n$ on $k^n$, inequality (\ref{grasstability}) holds.
\end{corollary}

These results imply that the preimages of $(\mathbb{P}_{m - 1})^{n, ss}(\GL_m)$ and
$\gr(m, n)^{ss}(T_n)$ in $\mathbb{P}\mathbb{M}_{n, m}$ coincide and we denote
this set by $\mathbb{P}\mathbb{M}_{n, m}^o$. Then the sets $(\mathbb{P}_{m - 1})^{n, ss}(\GL_m)$
and $\gr(m,$ $n)^{ss}(T_n)$ both are geometric quotients of
$\mathbb{P}\mathbb{M}_{n, m}^o$ and their quotients
\begin{equation*}
\gr(m, n)^{ss}(T_n) // T_n \text{\quad and \quad} (\mathbb{P}_{m - 1})^{n,ss}(\GL_m) // \GL_m.
\end{equation*}
are good quotients of $\mathbb{P}\mathbb{M}_{n, m}^o$ as each is a good quotient of a good
quotient.
By the universal property of good quotients, these two spaces coincide with the good
quotient $\mathbb{P}\mathbb{M}^o_{n,m} // G$, and thus are isomorphic.
In particular, there is a commutative diagram consisting of good quotients:
\begin{equation*}
\xymatrix{
& \mathbb{P}\mathbb{M}_{n, m}^o \ar[ld] \ar[rd] \\
(\mathbb{P}_{m - 1})^{n, ss}(\GL_m) \ar[rd] & & \gr(m, n)^{ss}(T_n) \ar[ld] \\
& \mathcal{M}_{n, m}
}
\end{equation*}
We want to extend this correspondence using the well-known isomorphism
\begin{equation*}
\gr(m, n) \cong \gr(n - m, n)
\end{equation*}
which can be interpreted as saying that an $n \times m$-matrix $A$ of rank $m$ representing
a point in $\gr(m, n)$ is mapped to an $(n - m) \times n$-matrix $\check{A}$ representing a
point in $\gr(n - m, m)$ such that both matrices fit into a short exact sequence
\begin{equation*}
0 \longrightarrow k^m \overset{A}{\longrightarrow} k^n
\overset{\check{A}}{\longrightarrow} k^{n - m} \longrightarrow 0
\end{equation*}
This correspondence is compatible with the action of the torus $T_n$ on both sides:
\begin{lemma}
Let $T_n$ be a maximal subtorus of $\GL_n$. Consider the actions of $T_n$ on $\gr(m, n)$
and $\gr(n - m, n)$, induced by its natural actions on $k^n$ and the dual vector space
$(k^n)\check{\ }$, respectively. Then the canonical isomorphism between $\gr(m, n)$ and
$\gr(n - m, n)$, which is induced by the canonical isomorphism between $k^n$ and
$(k^n)\check{\ }$, is $T_n$-equivariant and maps the (semi-)stable points as specified in
Proposition \ref{tuplequot}, to (semi-)stable points. There exists a natural isomorphism
\begin{equation*}
\gr(m, n)^{ss}(T_n) // T_n \cong \gr(n - m, n)^{ss}(T_n) // T_n 
\end{equation*}
\end{lemma}

\begin{proof}
A little bit of linear algebra shows that for some $A \in \gr(m, n)$ and for all
linear subspaces $L \subset k^n$ the following holds,
\begin{equation*}
\dim A \cap L < \frac{m}{n} \dim L \text{\quad if and only if \quad}
\dim A\check{\ } \cap L\check{\ } < \frac{n - m}{n} \dim L\check{\ }
\end{equation*}
(respectively $\leq$), where $A\check{\ }$ and $L\check{\ }$ are the annihilators
of $A$ and $L$ in $(k^n)\check{\ }$.
\end{proof}

Because of this we can extend our correspondences to the following diagram:
\begin{equation*}
\xymatrix{
& \mathbb{P}\mathbb{M}_{n, m}^o \ar[ld] \ar[rd] & & \mathbb{P}\mathbb{M}_{n - m, n}^o
\ar[ld], \ar[rd] \\
(\mathbb{P}_{m - 1})^{n,ss} \ar[rd] & & \gr(m, n)^{ss} \ar[ld] \cong
\gr(n - m, n)^{ss} \ar[rd] & &
(\mathbb{P}_{n - m - 1})^{n,ss} \ar[ld] \\
& \mathcal{M}_{n, m}  & \cong \quad \ \ & \mathcal{M}_{n - m, n}
}
\end{equation*}

\section{Moduli of Equivariant Sheaves}
\label{moduli}

Let us fix a tuple of nonnegative numbers $I = (i^\rho \mid \rho \in \rays)$ and a partition
$\mathcal{P} = \{\Pi_1, \dots, \Pi_s\}$ of the set $\Pi = \{\rho \in \rays \mid i^\rho > 0\}$ where
$s \geq 2$. In Section \ref{resolutions} we have identified such data as a set of typical discrete
parameters for  equivariant vector bundles of rank 2 on a toric surface $X$.

We have shown in Theorem \ref{evenmoreeulersequences} that for each such bundle \msh{E}
whose equivariant first Chern class in \mweildivisors\ and whose coarse partition
$\mathcal{P}_\mathcal{E}$ coincide with $I$ and $\mathcal{P}$, respectively,
there exists a short exact sequence of the form
\begin{equation*}
0 \longrightarrow \sh{O}^{s - 2} \overset{A}{\longrightarrow}
\bigoplus_{i = 1}^s \sh{O}(\sum_{\rho \in \Pi_i} i^\rho \cdot D_\rho)
\overset{\check{A}}{\longrightarrow} \sh{E} \longrightarrow 0
\end{equation*}
which corresponds to a short exact sequence of vector spaces
\begin{equation*}
0 \longrightarrow k^{s - 2} \overset{A^0}{\longrightarrow} k^s
\overset{\check{A}^0}{\longrightarrow} k^2 \longrightarrow 0.
\end{equation*}
In order to obtain moduli spaces, we ask for spaces which parametrize
isomorphism classes of equivariant vector bundles \msh{E} of rank 2 with fixed
coarse partition $\mathcal{P}_\mathcal{E} = \mathcal{P}$.
The conditions on $A$ in Theorem \ref{evenmoreeulersequences} imply that the set of
matrices $A$ whose cokernel is such a vector bundle is dense in $\mathbb{M}_{s, s - 2}$.
So by varying matrices $A$ we have a
natural candidate for a parameter space of vector bundles \msh{E} with fixed $I$ and
$\mathcal{P}$ which is given by $\mathbb{M}_{s, s - 2}$ modulo the equivariant automorphisms
of $\sh{O}^{s - 2}$ and $\bigoplus_{i = 1}^s \sh{O}(\sum_{\rho \in \Pi_i} i^\rho \cdot D_\rho)$,
$\GL_{s - 2}$ and $T_s$, respectively.

Another natural parameter space is the set of configurations of $s$ points in
$\mathbb{P}\mathbf{E}^0 \cong \mathbb{P}_1$ which can be given by the columns of the
matrix $\check{A}^0$. In that case equivariant isomorphism
classes of bundles are determined by configurations modulo linear transformations
by $\GL_2$. This is the sort of moduli space which has already been
suggested by Klyachko in \cite{Kly90}.

By the results of the previous section both spaces can be compared
in terms of the GIT-quotients $\mathcal{M}_{s, s - 2}$ and $\mathcal{M}_{2, s}$.
By Theorem \ref{desing} the isomorphism $\mathcal{M}_{s, s - 2} \overset{\cong}{\longrightarrow}
\mathcal{M}_{2, s}$ which is given by the map $A \mapsto \check{A}$, respectively
$A^0 \mapsto \check{A}^0$, can be interpreted as the map
\begin{equation*}
\sh{E} \mapsto \sh{E}\check{\ }\check{\ }
\end{equation*}
which is defined for any cokernel \msh{E} represented by some GIT-semistable matrix $A$ in
$\mathbb{P}\mathbb{M}_{s, s - 2}$.

Let us now investigate the semistable points of
$\mathbb{P}\mathbb{M}_{s, s - 2}$ and $\mathbb{P}\mathbb{M}_{2, s}$.
Recall from proposition \ref{tuplequot} that a point $(p_1, \dots, p_s) \in
(\mathbb{P}_1)^s$ is properly semistable with respect to the action of $\GL_2$
iff precisely $\frac{s}{2}$ of the $p_i$ coincide. Thus properly semistable points
exist only in the case $s = 2t$ even.

\begin{proposition}
\label{equivext}
Let $\sh{E} = \sh{E}(I, \mathcal{P}, A)$ be a torsion free sheaf given by a short
exact sequence as above.
Let $A_{i_1}$ the $i_1$-th column of $\check{A}$, let $\mathcal{P}_1 = \{\Pi_{i_1}, \dots
\Pi_{i_r}\}$ be the maximal subset of $\mathcal{P}$ with $\langle A_{i_k} \rangle = \langle
A_{i_1} \rangle$ for $1 \leq k \leq r$, and let
$\mathcal{P}_2$ be the complement of $\mathcal{P}_1$ in $\mathcal{P}$.
Then the torsion free sheaf $\sh{E}$ defined by $A$ is an extension
\begin{equation*}
0 \longrightarrow \sh{E}_1 \longrightarrow \sh{E} \longrightarrow \sh{E}_2 \longrightarrow 0
\end{equation*}
where $\sh{E}_1$ and $\sh{E}_2$ are torsion free sheaves of rank 1 with
\begin{equation*}
\sh{E}_1\check{\ }\check{\ } \cong \sh{O}(\sum_{\Pi \in \mathcal{P}_1}
\sum_{\rho \in \Pi}i^\rho \cdot D_\rho)
\text{\quad and \quad}
\sh{E}_2\check{\ }\check{\ } \cong \sh{O}(\sum_{\Pi \in \mathcal{P}_2}
\sum_{\rho \in \Pi}i^\rho \cdot D_\rho).
\end{equation*}
\end{proposition}

\begin{proof}
We obtain this extension by partition of the matrix $\check{A}$ via the following
diagram:
\begin{equation*}
\xymatrix{
& 0 \ar[d] & 0 \ar[d] & 0 \ar[d] \\
0 \ar[r] & \sh{O}^{r - 1} \ar[r]^-{A_1} \ar[d] &  \bigoplus_{\Pi \in \mathcal{P}_1}
\sh{O}(\sum_{\rho \in \Pi} i^\rho \cdot D_\rho) \ar[r]^-{\check{A}_1} \ar[d] & \sh{E}_1
\ar[r] \ar[d] & 0 \\
0 \ar[r] & \sh{O}^{s - 2} \ar[r]^-A \ar[d] & \bigoplus_{\Pi \in \mathcal{P}}
\sh{O}(\sum_{\rho \in \Pi} i^\rho \cdot D_\rho) \ar[r]^-{\check{A}} \ar[d] & \sh{E} \ar[r] \ar[d] & 0 \\
0 \ar[r] & \sh{O}^{s - r - 1} \ar[r]^-{A_2} \ar[d] & \bigoplus_{\Pi \in \mathcal{P}_2}
\sh{O}(\sum_{\rho \in \Pi} i^\rho \cdot D_\rho) \ar[r]^-{\check{A}_2} \ar[d] & \sh{E}_2 \ar[r] \ar[d] & 0 \\
& 0 & 0 & 0
}
\end{equation*}
where $A_1$ is represented by the submatrix of $A$ consisting of the rows corresponding to
$\mathcal{P}_1$.
\end{proof}

\begin{corollary}
\label{equivextcor}
Let \msh{E} and $\mathcal{P}$ be as above and let $\sh{F} \subset \sh{E}$ be any torsion free
equivariant subsheaf of rank 1. Then there exists a subset $\mathcal{P}' \subset
\mathcal{P}$ such that $\sh{F}\check{\ }\check{\ } \cong \sh{O}(\sum_{\Pi \in \mathcal{P}'}\sum_{\rho
\in \Pi} i^\rho \cdot D_\rho)$.
\end{corollary}

\begin{corollary}
\label{extequiv}
Let $s = 2t$ be even, let $\check{A}$ represent a properly semistable point in
$(\mathbb{P}_1)^{ss}$ and let $0 \longrightarrow \sh{E}_1 \longrightarrow \sh{E}
\longrightarrow \sh{E}_2 \longrightarrow 0$ be the corresponding extension.
Then the image of $A$ in $\mathcal{M}_{s, s - 2}$ represents all matrices whose
corresponding extensions are in $\operatorname{Ext}^1(\sh{E}_2, \sh{E}_1)$ or
$\operatorname{Ext}^1(\sh{E}_1, \sh{E}_2)$, i.e. \msh{E} is GIT-equivalent to
the direct sum $\sh{E}_1 \oplus \sh{E}_2$.
\end{corollary}

\begin{proof}
Each orbit in $(\mathbb{P}_1)^s$ which contains a point $(p_1, \dots, p_s)$ such that
some $t$ points $p_{i_1}, \dots p_{i_t}$ coincide contains in its closure the points
of the form $(p_{i_1}, \dots,$ $p_{i_t}, \dots p, \dots p)$ for some $p_{i_1} \neq p \in
\mathbb{P}_1$.
\end{proof}

In the generic case, we have in particular:

\begin{corollary}
Let $n = \# \rays$ be even and $i^\rho > 0$ for all $\rho \in \rays$, let
$\mathcal{P} = \{\{\rho\}\}_{\rho \in \rays}$ be the fine partition of \mrays,
and let $n = \# \mathcal{P}$. Then there exists precisely
one point in $\mathcal{M}_{n, n - 2}$ which can be represented by a direct sum
$\sh{E}_1 \oplus \sh{E}_2$ such that $\sh{E}_1$ and $\sh{E}_2$ are locally free.
\end{corollary}

\begin{proof}
There exists precisely one partition $\Pi_1 \dot{\cup} \Pi_2 = \rays$ such that the
$\Pi_i$ do not contain two adjacent rays and which is given by
$\Pi_1 = \{\{\rho_{2 \cdot i}\} \vert 1 \leq i \leq \frac{n}{2}\}$.
\end{proof}

The observation made in Proposition \ref{equivext} motivates the following definition:
\begin{definition}
\label{stabledef}
Let $I$ and $\mathcal{P}$ as before and let \msh{E} be a torsion free equivariant sheaf of rank
2 over $X$ such that $\mathcal{P}$ is a refinement of the coarse partition
$\mathcal{P}_{\sh{E}\check{\ }\check{ \ }}$ associated to
the locally free sheaf $\sh{E}\check{\ }\check{\ }$.
Let $\sh{F} \subset \sh{E}$ be a torsion free equivariant subsheaf of rank 1. Then by
\ref{equivextcor}
$\sh{F}\check{\ }\check{\ } \cong \sh{O}(\sum_{\Pi \in \mathcal{P}'}\sum_{\rho \in \Pi} i^\rho \cdot D_\rho)$ with a unique subset $\mathcal{P}' \subset \mathcal{P}$.
We say that \msh{E} is {\it $\mathcal{P}$-stable} (respectively {\it $\mathcal{P}$-semistable}) if for every
equivariant torsion free subsheaf $\sh{F} \subset \sh{E}$ of rank 1
$\# \mathcal{P}' < \frac{1}{2} \# \mathcal{P}$ (respectively $\# \mathcal{P}' \leq \frac{1}{2} \# \mathcal{P}$).
\end{definition}

\begin{theorem}
Let $i_\rho > 0$ for $\rho \in \Pi \subset \rays$ and let $\mathcal{P} = \{\Pi_1, \dots,
\Pi_s\}$ be a partition of $\Pi$. Consider short exact sequences
\begin{equation*}
0 \longrightarrow \sh{O}^{s - 2} \overset{A}{\longrightarrow}
\bigoplus_{i = 1}^s \sh{O}(\sum_{\rho \in \Pi_i} i^\rho \cdot D_\rho)
\longrightarrow \sh{E} \longrightarrow 0.
\end{equation*}
Then \msh{E} is {\it $\mathcal{P}$-stable} (respectively {\it $\mathcal{P}$-semistable}) if and only
if $A$ represents a GIT-stable (respectively GIT-semistable) point in $(\mathbb{P}_{s - 3})^s$
with respect to the action of $\GL_{s - 2}$.
\end{theorem}

\begin{proof}
This follows from the fact that we can represent $A$ by a configuration $(p_1, \dots, p_s)$ of
points in $(\mathbb{P}_1)^s$. Then Definition \ref{stabledef} is equivalent to the fact that at most
$\frac{s}{2}$ of the points $p_i$ coincide.
\end{proof}

Now we can define an equivalence relation on the set of
$\mathcal{P}$-semistable sheaves as follows:

\begin{definition}
Let \msh{E} and $\sh{E}'$ be $\mathcal{P}$-semistable sheaves. Then we say that \msh{E}
and $\sh{E}'$ are $\mathcal{P}$-equivalent iff one of the following conditions holds:
\begin{enumerate}[(i)]
\item \msh{E} and $\sh{E}'$ both are $\mathcal{P}$-stable and equivariantly isomorphic,
$\sh{E} \cong \sh{E}'$,
\item  \msh{E} and $\sh{E}'$ both are $\mathcal{P}$-semistable and the following holds.
Let $\Psi \in \mathcal{P}_\mathcal{E}$ such that $\# \Psi = \frac{s}{2}$. Then either $\Psi
\in \mathcal{P}_{\mathcal{E}'}$ or $\Pi \setminus \Psi \in \mathcal{P}_{\mathcal{E}'}$.
\end{enumerate}
\end{definition}

The last condition implies that if $0 \longrightarrow \sh{E}_1 \longrightarrow \sh{E}
\longrightarrow \sh{E}_2 \longrightarrow 0$ is the extension of $\sh{E}$ corresponding
to $\Psi$ as in Proposition \ref{equivext}, then by Corollary \ref{extequiv} $\sh{E}$
is $\mathcal{P}$-equivalent to $\sh{E}_1 \oplus \sh{E}_2$. From this definition follows

\begin{theorem}
Fix numbers $\{i^\rho \geq 0\}_{\rho \in \rays}$, let $\Pi = \{\rho \mid i^\rho > 0\} \subset \rays$
and let $\mathcal{P} = \{\Pi_1, \dots, \Pi_s\}$. Then $\mathcal{M}_{s, s - 2}$ is the set of
$\mathcal{P}$-equivalence classes of $\mathcal{P}$-semistable torsion free equivariant
sheaves of rank 2 on $X$.
\end{theorem}

\begin{definition}
If $\mathcal{P}$ is fixed, we denote $\mathcal{M}_\mathcal{P} := \mathcal{M}_{s, s - 2}$ and
call $\mathcal{M}_\mathcal{P}$ {\it moduli space of $\mathcal{P}$-equivalence classes}.
\end{definition}

It remains to show that the spaces $\mathcal{M}_\mathcal{P}$ are moduli spaces of
suitably defined $\mathcal{P}$-families, e.g. in the sense of \cite{Newstead}. A detailed
treatment of this problem would require to generalize all our constructions to families.
We hope to come back to this question in a more general context in future work.

\begin{remark}
There is the following result of Klyachko:

\begin{proposition}[\cite{Kly90}, Corollary 1.2.5]
Let \msh{E} and $\sh{E}'$ be two equivariant vector bundles over a smooth complete toric
variety. If there exists an arbitrary isomorphism $\sh{E} \cong \sh{E}'$ of vector
bundles, then there is an $m \in M$ such that there is an equivariant
isomorphism $\sh{E} \cong \sh{E} \otimes \sh{O}(\chi(m))$, where $\sh{O}(\chi(m))$ denotes
the structure sheaf endowed with the action by the character $\chi(m)$.
\end{proposition}

This means that in our situation, where $X$ is complete, equivariant isomorphism classes
and isomorphism classes of vector bundles coincide up to a twist with a character.
So after fixing numbers $i^\rho$, $\rho \in \rays$, the subspace $\mathcal{M}'_\mathcal{P}
\subset \mathcal{M}_\mathcal{P}$ consisting of isomorphism classes of vector bundles
even classifies non-equivariant isomorphism classes.
\end{remark}

\begin{remark}
We want to point out that our moduli depend only on the {\it combinatorial} structure
of the underlying toric variety, that is, the number of rays \mrays\ in the fan of $X$,
but not on the concrete realization of the fan $\Delta$ inside the lattice $N$.
\end{remark}

\begin{example}
Let $X = \mathbb{P}_2(k)$, then we have to consider the quotient of
$\mathbb{P}\mathbb{M}_{3, 1}$
by the group $G \cong \GL_1 \times T_3 \cong k^* \times (k^*)^3$. This quotient is just
a point, i.e. the set of equivariant isomorphism classes of indecomposable equivariant
vector bundles of rank 2 on $\mathbb{P}_2(k)$ is {\it discrete}. This reproduces the
original result of Kaneyama (\cite{Kan1}).
\end{example}

\begin{example}
Let $a \geq 0$ and $X = \mathbb{F}_a$ a Hirzebruch surface. Assume that
the rays $\rho_1, \dots, \rho_4$ are enumerated clockwise. The set
$(\mathbb{P}_1)^{4, s}$ of stable points of $(\mathbb{P}_1)^{4}$ with respect
to the diagonal action of $\GL_2$ is
\begin{equation*}
\{(p_1, \dots, p_4) \subset (\mathbb{P}_1)^{4} \mid p_i \neq p_j \text{ for all } i \neq j\},
\end{equation*}
i.e. the set of four-point configurations in $\mathbb{P}_1$ no two points of which coincide.
There is an isomorphism $(\mathbb{P}_1)^{4, s} \cong \PGL_2 \times (\mathbb{P}_1 \setminus
\{0, 1, \infty\})$ which we choose to be
\begin{equation*}
(p_1, p_2, p_3, p_4) \mapsto (g, g . p_4),
\end{equation*}
where $g \in \PGL_2$ is the unique element which moves the points $p_1$, $p_2$, $p_3$
to the positions $0$, $1$, and $\infty$, respectively. The inverse map is given by
\begin{equation*}
(g, p) \mapsto (g^{-1} 0, g^{-1} 1, g^{-1} \infty, g^{-1} p).
\end{equation*}
The quotient $\mathbb{P}\mathbb{M}_{4, 2}^{ss}(\GL_2 \times T_2) // \GL_2 \times T_2$ has
a completion by semistable points:
\begin{equation*}
(\mathbb{P}_1)^{4, ss} = \{ (p_1, p_2, p_3, p_4) | \text{ such that no three points
$p_i$ coincide}\}
\end{equation*}

In terms of $4 \times 2$ matrices this means that each semistable but not stable matrix
can be brought into one of six standard forms with at most one zero in a row and two zeros
in a column:
\begin{equation*}
A_{34} = \begin{pmatrix}
* & * \\
* & * \\
0 & * \\
0 & * 
\end{pmatrix}, \quad
A_{24} = \begin{pmatrix}
* & * \\
0 & * \\
* & * \\
0 & *
\end{pmatrix},
\quad \dots \quad,
A_{12} = \begin{pmatrix}
0 & * \\
0 & * \\
* & * \\
* & *
\end{pmatrix}.
\end{equation*}
The image of a matrix $A_{ij}$ in $\mathcal{M}_{4, 2}$ represents the $\mathcal{P}$-equivalence
class of $\sh{E}_1 \oplus \sh{E}_2$ where $\sh{E}_1\check{\ }\check{\ }
\cong \sh{O}(i^{\rho_i} \cdot D_{\rho_i} + i^{\rho_j} \cdot D_{\rho_j})$
and $\sh{E}_2\check{\ }\check{\ } \cong \sh{O}(i^{\rho_k} \cdot D_{\rho_k} +
i^{\rho_l} \cdot D_{\rho_l})$ with $\{i, j, k, l\} = \{0, 1, 2, 3\}$.
By our choice of coordinates the matrices of types $A_{13}$ and $A_{24}$ with locally
free cokernels are mapped to the point $1 \in \mathbb{P}_1$.
\end{example}

\begin{remark}
Let $s > 4$, let $X_{ij}$, $1 \leq i \leq s$, $1 \leq j \leq s - 2$ be the coordinates of
$\mathbb{M}_{s, s - 2}$ and let $\mathcal{X} = (X_{ij})$. Then the determinants of
the minors $\mathcal{X}^{ij}$ of
$\mathcal{X}$ describe $\frac{1}{2} s (s - 1)$ hypersurfaces $T^{ij}_{\mathcal{P}}$ in
$\mathbb{M}_{s, s - 2}$. Via the duality of Section \ref{duality}, the minors
$\mathcal{X}^{ij}$ describe precisely the configurations of points $(p_1, \dots, p_s)$
in $\mathbb{P}_1$,
where the points $p_i$ and $p_j$ coincide. Denote $\mathcal{T}_{\mathcal{P}}^{ij}$ the
image of $T^{ij}_{\mathcal{P}}$ in $\mathcal{M}_\mathcal{P}$. $\mathcal{M}_\mathcal{P}$
is a good quotient of $\mathbb{P}\mathbb{M}_{s, s - 2}$, and thus
$\mathcal{T}_{\mathcal{P}}^{ij}$
is a closed subset of $\mathcal{M}_\mathcal{P}$. Moreover, because $s > 4$, by Proposition
\ref{tuplequot} we see that each $\mathcal{T}_{\mathcal{P}}^{ij}$ contains a dense
subset whose preimage in $\mathbb{P}\mathbb{M}_{s, s - 2}$ consists of stable points. Thus,
the $\mathcal{T}^{ij}_{\mathcal{P}}$ describe $\frac{1}{2} s (s - 1)$
hypersurfaces of $\mathcal{M}_\mathcal{P}$ and the locus in $\mathcal{M}_\mathcal{P}$
consisting of torsion free sheaves which are not locally free is described by the
$s$ divisors $\mathcal{T}^{s,1}$ and $\mathcal{T}^{i, i + 1}$ for $1 \leq i \leq s - 1$.
\end{remark}

\begin{example}
Let $X$ be a toric surface which has six rays. The space $\mathcal{M}_{2,6}$
has been calculated in \cite{Dolgachev} and is a cubic hypersurface in
$\mathbb{P}_4$ defined by the equation
$X_1X_2X_4-X_3X_0X_4+X_3X_1X_2+X_3X_0X_1+X_3X_0X_2-X_3X_0^2 = 0$.
This hypersurface has ten nodes representing precisely the ten $\mathcal{P}$-equivalence
classes of $\mathcal{P}$-semistable but not $\mathcal{P}$-stable sheaves.
\end{example}

\end{document}